\documentclass{gtart}

\def\ifplaintex{\expandafter\ifx\csname documentclass\endcsname\relax}

\ifplaintex 
\else
\fi

\expandafter\ifx\csname beginpicture\endcsname\relax
\expandafter\ifx\csname documentclass\endcsname\relax
\input pictex \else
\input prepictex \input pictex \input postpictex \fi\fi

\def\gt{{\mathsurround=0pt\it $\cal G\mskip-2mu$eometry \&\ 
$\cal T\!\!$opology}}        

\def\gtp{{\mathsurround=0pt\it $\cal G\mskip-2mu$eometry \&\ 
$\cal T\!\!$opology $\cal P\!$ublications}}  

\def\lognumber#1{\def\thelognumber{#1}}
\def\volumenumber#1{\def\thevolumenumber{#1}}
\def\papernumber#1{\def\thepapernumber{#1}}
\def\volumeyear#1{\def\thevolumeyear{#1}}

\def\pagenumbers#1#2{\def\startpage{#1}\def\finishpage{#2}}
\def\published#1{\def\publishdate{#1}}
\def\proposed#1{\def\theproposer{#1}}
\def\seconded#1{\def\theseconders{#1}}
\def\received#1{\def\receiveddate{#1}}

\def\accepted#1{\def\accepteddate{#1}}
\def\asciititle#1{\def\theasciititle{#1}}

\def\asciiaddress#1{\def\theasciiaddress{#1}}

\long\def\asciiabstract#1{\long\def\theasciiabstract{#1}}
\def\asciikeywords#1{\def\theasciikeywords{#1}}

\let\\\par\let\thelognumber\relax
\let\thevolumenumber\relax\let\thepapernumber\relax
\let\thevolumeyear\relax\let\thesamplenumber\relax\let\startpage\relax
\let\finishpage\relax\let\publishdate\relax\let\receiveddate\relax
\let\reviseddate\relax\let\accepteddate\relax\let\theasciititle\relax
\let\theasciiauthors\relax\let\theasciiaddress\relax
\let\theasciiabstract\relax\let\theasciikeywords\relax
\let\theasciiemail\relax\let\theshortauthors\relax\let\theshorttitle\relax

\long\def\maketitlep{   

\count0=\startpage

\gt\hfill      
\beginpicture
\setcoordinatesystem units <0.33truein, 0.33truein> point at 2.2 0.9
\setplotsymbol ({$\cal G$})
\plotsymbolspacing=9truept
\circulararc 315 degrees from 0 1 center at 0 0
\setplotsymbol ({$\cal T$})
\circulararc 315 degrees from 1 -1 center at 1 0
\endpicture
\break
{\small\ifx\thesamplenumber\relax 
Volume \else Sample
\fi\thevolumenumber\ (\thevolumeyear)
\startpage--\finishpage\nl
Published: \publishdate}
\vglue 0.5truein plus 0.4fil minus 0.1truein

{\parskip=0pt\leftskip 0pt plus 1fil\def\\{\par\smallskip}{\ifplaintex\large
\else\Large\fi\bf\thetitle}\par\medskip}   

\vglue 0pt plus 0.1fil 

{\parskip=0pt\leftskip 0pt plus 1fil\def\\{\par}{\sc\theauthors}
\par\medskip}

\vglue 0pt plus 0.1fil 

{\small\parskip=0pt\let\newline\\
{\leftskip 0pt plus 1fil\def\\{\par}{\sl\theaddress}\par}
\expandafter\ifx\theemail\relax    
\relax\else\vglue 5pt plus 0.02fil minus 2pt\def\\{\stdspace{\rm 
and}\stdspace} 
\cl{Email:\stdspace\tt\theemail}\fi
\ifx\theurl\relax                  
\relax\else\vglue 5pt plus 0.02fil minus 2pt\def\\{\stdspace{\rm 
and}\stdspace}
\cl{URL:\stdspace\tt\theurl}\fi\par}

\vglue 7pt plus 0.3fil minus 3pt

{\bf Abstract}
\vglue 5pt plus 0.1fil minus 2pt

\theabstract

\vglue 7pt plus 0.3fil minus 3pt

{\bf AMS Classification numbers}\quad Primary:\quad \theprimaryclass

Secondary:\quad \thesecondaryclass

\vglue 5pt plus 0.3fil minus 2pt

{\bf Keywords}\quad \thekeywords

\vglue 10pt plus 0.5fil minus 5pt

{\small  Proposed: \theproposer\hfill Received: \receiveddate\nl
Seconded: \theseconders\hfill 
\ifx\reviseddate\relax                         
Accepted: \accepteddate                        
\else
Revised: \reviseddate                          
\fi}
\eject
}       

\let\maketitlepage\maketitlep
\let\maketitle\maketitlepage

\font\phead=cmsl9 scaled 950
\font=cmsl9 scaled 1050
\font\pnum=cmbx10 scaled 913
\font\lnum=cmbx10 
\font\pfoot=cmsl9 scaled 950
\font=cmsl9 scaled 1050
\ifplaintex
\headline{\vbox to 0pt{\vskip -4.5mm\line{\small\phead\ifnum
\count0=\startpage ISSN 1364-0380 (on line)
1465-3060 (printed) \hfill {\pnum\folio}\else\ifodd\count0\def\\{ }%
\ifx\theshorttitle\relax\thetitle\else\theshorttitle\fi\hfill{\pnum\folio}
\else\def\\{ and }{\pnum\folio}\hfill\ifx\theshortauthors\relax\theauthors
\else\theshortauthors\fi\fi\fi}\vss}}
\footline{\vbox to 0pt{\vglue 0mm\line{\small\pfoot\ifnum\count0=\startpage
\copyright\ \gtp\hfill\else
\gt, Volume \thevolumenumber\ (\thevolumeyear)\hfill\fi}\vss
}}
\else
\makeatletter
\def\@oddhead{{\small\lhead\ifnum\count0=\startpage ISSN 1364-0380 (on line)
1465-3060 (printed) \hfill {\lnum\number\count0}\else\ifodd\count0
\def\\{ }\ifx\theshorttitle\relax \thetitle \else\theshorttitle\fi\hfill
{\lnum\number\count0}\else\def\\{ and }{\lnum\number\count0}
\hfill\ifx\theshortauthors\relax 
\theauthors\else\theshortauthors\fi\fi\fi}}\def\@evenhead{\@oddhead}
\def\@oddfoot{\small\lfoot\ifnum\count0=\startpage\copyright\ \gtp\hfill\else
\gt, Volume \thevolumenumber\ (\thevolumeyear)\hfill\fi}
\def\@evenfoot{\@oddfoot}
\makeatother
\fi

\newwrite\gtoutfile
\long\gdef\makeheadfile{  
{\def\\{, }\def\s{ }
\immediate\openout\gtoutfile head.xxx
\immediate\write\gtoutfile{Proxy-for: \ifx\theasciiauthors\relax
\theauthors\else\theasciiauthors\fi\s<\ifx\theasciiemail\relax\theemail\else\theasciiemail\fi>}
\immediate\write\gtoutfile{\noexpand\\}
\immediate\write\gtoutfile{Authors: \ifx\theasciiauthors\relax
\theauthors\else\theasciiauthors\fi}
{\def\\{ }\immediate\write\gtoutfile{Title: \ifx\theasciititle\relax
\thetitle\else\theasciititle\fi}}
\immediate\write\gtoutfile{Subj-class: GT or SG or MG etc}
\immediate\write\gtoutfile{MSC-class: \theprimaryclass\ifx\thesecondaryclass\relax\else, \thesecondaryclass\fi}
\immediate\write\gtoutfile{Journal-ref: Geom. Topol. \thevolumenumber
(\thevolumeyear) \startpage-\finishpage}
\immediate\write\gtoutfile{Comments: Published by Geometry and Topology at}
\immediate\write\gtoutfile{\s\s http://www.maths.warwick.ac.uk/gt/GTVol\thevolumenumber/paper\thepapernumber.abs.html}
\immediate\write\gtoutfile{\noexpand\\}
\immediate\write\gtoutfile{}
\ifx\theasciiabstract\relax
\immediate\write\gtoutfile{\theabstract}\else
\immediate\write\gtoutfile{\theasciiabstract}\fi
\immediate\write\gtoutfile{}
\immediate\write\gtoutfile{\noexpand\\}
\immediate\write\gtoutfile{}
\immediate\closeout\gtoutfile}}  

\def\maketitlepage{\maketitlep\makeheadfile}
\let\maketitle\maketitlepage

\lognumber{327}
\volumenumber{7}\papernumber{24}\volumeyear{2003}
\pagenumbers{789}{797}
\received{3 June 2003}
\published{28 November 2003}
\accepted{13 November 2003}

\proposed{Jean-Pierre Otal}
\seconded{David Gabai, Benson Farb}

\usepackage{amsmath,amssymb}

\newtheorem{sat}{Theorem}
\newtheorem*{sat*}{Theorem}
\newtheorem{lem}[sat]{Lemma}

\newtheorem{prop}[sat]{Proposition}

\newtheorem*{namedtheorem}{\theoremname}
\newcommand{\theoremname}{testing}
\newenvironment{named}[1]{\renewcommand{\theoremname}{#1}\begin{namedtheorem}}{\end{namedtheorem}}

\newcommand{\D}{\partial}
\newcommand{\defin}{\ensuremath{\overset{ \text{\tiny def} }{=} } }
\newcommand{\BC}{\mathbb C}
\newcommand{\BH}{\mathbb H}
\newcommand{\BR}{\mathbb R}

\newcommand{\BS}{\mathbb S}
\newcommand{\BZ}{\mathbb Z}
\newcommand{\CA}{\mathcal A}
\newcommand{\CL}{\mathcal L}
\newcommand{\CM}{\mathcal M}
\newcommand{\CN}{\mathcal N}
\newcommand{\CO}{\mathcal O}
\newcommand{\CP}{\mathcal P}
\newcommand{\CU}{\mathcal U}

\DeclareMathOperator{\supp}{supp}

\DeclareMathOperator{\inter}{Int}
\DeclareMathOperator{\hypvol}{Hypvol}
\DeclareMathOperator{\hyp}{hyp}
\DeclareMathOperator{\Minvol}{Minvol}
\DeclareMathOperator{\vol}{vol}

\begin{document}
\title{Hyperbolic cone--manifolds with large cone--angles}
\asciititle{Hyperbolic cone-manifolds with large cone-angles}
\author{Juan Souto}
\address{Mathematisches Institut\\Universit\"at Bonn\\Beringstr. 1, 53115 Bonn, Germany}
\asciiaddress{Mathematisches Institut, Universitaet 
Bonn\\Beringstr. 1, 53115 Bonn, Germany}
\email{souto@math.uni-bonn.de}
\url{http://www.math.uni-bonn.de/people/souto}
\begin{abstract}
We prove that every closed oriented 3--manifold admits a hyperbolic cone--manifold structure with cone--angle arbitrarily close to $2\pi$.
\end{abstract}
\asciiabstract{%
We prove that every closed oriented 3-manifold admits a hyperbolic
cone-manifold structure with cone-angle arbitrarily close to 2pi.}

\asciikeywords{Hyperbolic cone-manifold, Kleinian groups}
\maketitle 

\section{Introduction}

Consider the hyperbolic 3--space in the upper half--space model $\BH^3\simeq\BC\times\BR_+$ and for $\alpha\in(0,2\pi)$ set $S_\alpha=\{(e^{r+i\theta},t)\mid\theta\in[0,\alpha], r\in[0,\infty), t\in\BR_+\}$. The boundary of $S_\alpha$ is a union of two hyperbolic half--planes. Denote by $\BH^3(\alpha)$ the space obtained from $S_\alpha$ by identifying both half--planes by a rotation around the vertical line $\{0\}\times\BR_+$. 

A distance on a 3--manifold $M$ determines a hyperbolic cone--manifold structure with singular locus a link $L\subset M$ and cone--angle $\alpha\in(0,2\pi)$, if every point $x\in M$ has a neighborhood which can be isometrically embedded either in $\BH^3$ or in $\BH^3(\alpha)$ depending on $x\in M\setminus L$ or $x\in L$.

Jean--Pierre Otal showed that the connected sum $\#^k(\BS^2\times\BS^1)$ of $k$ copies of $\BS^2\times\BS^1$ admits a hyperbolic cone--manifold structure with cone--angle $2\pi-\epsilon$ for all $\epsilon>0$ as follows: The manifold $\#^k(\BS^2\times\BS^1)$ is the double of the genus $k$ handlebody $H$. There is a convex--cocompact hyperbolic metric on the interior of $H$ such that the boundary of the convex--core is bent along a simple closed curve $\gamma$ with dihedral angle $\pi-\frac 12\epsilon$ \cite{Bonahon-Otal}; the convex--core is homeomorphic to $H$ and hence the double of the convex--core is homeomorphic to $\#^k(\BS^2\times\BS^1)$. The induced distance determines a hyperbolic cone--manifold structure on $\#^k(\BS^1\times\BS^2)$ with singular locus $\gamma$ and cone--angle $2\pi-\epsilon$. The same argument applies for every manifold which is the double of a compact manifold  whose interior admits a convex--cocompact hyperbolic metric. Michel Boileau asked whether every $3$--manifold has this property. Our goal is to give a positive answer to this question. We prove:

\begin{sat}\label{main}
Let $M$ be a closed and orientable $3$--manifold. For every $\epsilon$ there is a distance $d_\epsilon$ which determines a hyperbolic cone--manifold structure on $M$ with cone--angle $2\pi-\epsilon$.
\end{sat}

Before going further, we remark that we do not claim that the singular locus is independent of $\epsilon$. 

We now sketch the proof of Theorem \ref{main}. First, we construct a compact manifold $M^0$, whose boundary consists of tori, and such that there is a sequence $(M^0_n)$ of 3--manifolds obtained from $M^0$ by Dehn filling such that $M^0_n$ is homeomorphic to $M$ for all $n$. The especial structure of $M^0$ permits us to show that the interior $\inter M^0$ of the manifold $M^0$ admits, for every $\epsilon>0$, a complete hyperbolic cone--manifold structure with cone--angle $2\pi-\epsilon$. Thus, it follows from the work of Hodgson and Kerckhoff \cite{Hodgson-Kerckhoff} that for $n_\epsilon$ sufficiently large there is a distance $d_{n_\epsilon}=d_\epsilon$ on the manifold $M^0_{n_\epsilon}=M$ which determines a hyperbolic cone--manifold structure with cone--angle $2\pi-\epsilon$.

Let $(\epsilon_i)$ be a non-increasing sequence of positive numbers tending to $0$. If the corresponding sequence $(n_{\epsilon_i})$ grows fast enough, then the pointed Gromov--Hausdorff limit of the sequence $(M,d_{\epsilon_i})$ of metric spaces is a complete, smooth, hyperbolic manifold $X$ with finite volume. Moreover, the volume of the $(M,d_{\epsilon_i})$ converges to the volume of $X$ when $i$ tends to $\infty$; in particular the volume of $(M,d_{\epsilon_i})$ is uniformly bounded.

\medskip

I would like to thank Michel Boileau for many useful suggestions and remarks which have clearly improved the paper.

The author has been supported by the Sonderforschungsbereich 611.

\section{Preliminaries}\label{preli}

\subsection{Dehn filling}
Let $N$ be a compact manifold whose boundary consists of tori $T_1,\dots,T_k$ and let $U_1,\dots,U_k$ be solid tori. For any collection $\{\phi_i\}_{i=1,\dots,k}$ of homeomorphisms $\phi_i:\D U_i\to T_i$ let $N_{\phi_1,\dots,\phi_k}$ be the manifolds obtained from $N$ by attaching the solid torus $U_i$ via $\phi_i$ to $T_i$ for $i=1,\dots,k$.

Suppose that for all $i$ we have a basis $(m_i,l_i)$ of $H_1(T_i;\BZ)$ and let $\mu_i$ be the meridian of the solid torus $U_i$. There are coprime integers $a_i,b_i$ with $\phi_i(\mu_i)=a_im_i+b_il_i$ in $H_1(T_i;\BZ)$ for all $i=1,\dots,k$. It is well known that the manifold $N_{\phi_1,\dots,\phi_k}$ depends only on the set $\{a_1m_1+b_1l_1,\dots,a_km_k+b_kl_k\}$ of homology classes. We denote this manifold by $N_{(a_1m_1+b_1l_1),\dots,(a_km_k+b_kl_k)}$ and say that it has been obtained from $N$ filling the curves $a_im_i+b_il_i$.

The following theorem, due to Hodgson and Kerckhoff \cite{Hodgson-Kerckhoff} (see also \cite{Bromberg-rigidity}), generalizes Thurston's Dehn filling theorem:

\begin{named}{Generalized Dehn filling theorem}
Let $N$ be a compact manifold whose boundary consists of tori $T_1,\dots,T_k$ and let $(m_i,l_i)$ be a basis of $H_1(T_i;\BZ)$ for $i=1,\dots,k$. Assume that the interior $\inter N$ of $N$ admits a complete finite volume hyperbolic cone--manifold structure with cone--angle $\alpha\le 2\pi$. Then there exists $C>0$ with the following property:

The manifold $N_{(a_1m_1+b_1l_1),\dots,(a_km_k+b_kl_k)}$ admits a hyperbolic cone--manifold\break structure with cone--angle $\alpha$ if $\vert a_i\vert+\vert b_i\vert\ge C$ for all $i=1,\dots,k$.
\end{named}

\subsection{Geometrically finite manifolds}
The {\em convex--core} of a complete hyperbolic manifold $N$ with finitely generated fundamental group is the smallest closed convex set $CC(N)$ such that the inclusion $CC(N)\hookrightarrow N$ is a homotopy equivalence. The convex--core $CC(N)$ has empty interior if an only if $N$ is Fuchsian; since we will not be interested in this case we assume from now on that the interior of the convex--core is not empty. We will only work with geometrically finite manifolds, i.e. the convex--core has finite volume. If $N$ is geometrically finite then it is homeomorphic to the interior of a compact manifold $\CN$ and the convex--core $CC(N)$ is homeomorphic to $\CN\setminus\CP$ where $\CP\subset\D\CN$ is the union of all toroidal components of $\D\CN$ and of a collection of disjoint, non--parallel, essential simple closed curves. The pair $(\CN,\CP)$ is said to be the pared manifold associated to $N$ and $\CP$ is its parabolic locus (\cite{Japaner}). 

A theorem of Thurston \cite{Thurston} states that the induced distance on the boundary $\D CC(N)$ of the convex--core $CC(N)$ is a complete smooth hyperbolic metric with finite volume. The boundary components are in general not smoothly embedded, they are pleated surfaces bent along the so--called bending lamination. We will only consider geometrically finite manifolds for which the bending lamination is a weighted curve $\epsilon\cdot\gamma$. Here $\gamma$ is the simple closed geodesic of $N$ along which $\D CC(N)$ is bent and $\pi-\epsilon$ is the dihedral angle. 

The following theorem, due to Bonahon and Otal, is an especial case of \cite[Th\'eor\`eme 1]{Bonahon-Otal}.

\begin{named}{Realization theorem}
Let $\CN$ be a compact 3--manifold with incompressible boundary whose interior $\inter\CN$ admits a complete hyperbolic metric with parabolic locus $\CP$. If $\gamma\subset\D\CN\setminus\CP$ is an essential simple closed curve such that $\D\CN\setminus(\gamma\cup\CP)$ is acylindrical then for every $\epsilon>0$ there is a unique geometrically finite hyperbolic metric on $\inter\CN$ with parabolic locus $\CP$ and bending lamination $\epsilon\cdot\gamma$.
\end{named}

We refer to \cite{CEG} and to \cite{Japaner} for more about the geometry of the convex--core of geometrically finite manifolds.


\section{Proof of Theorem \ref{main}}\label{sec:main}

Let $S\subset M$ be a closed embedded surface which determines a Heegaard splitting $M=H_1\cup_\phi H_2$ of $M$. Here $H_1$ and $H_2$ are handlebodies and $\phi:\D H_1\to\D H_2$ is the attaching homeomorphism. Without loss of generality we may assume that $S$ has genus $g\ge 2$.

\begin{lem}\label{lem:OO}
There is a pant decomposition $P$ of $\D H_1$ such that both pared manifolds $(H_1,P)$ and $(H_2,\phi(P))$ have incompressible and acylindrical boundary.
\end{lem}
\begin{proof}
The Masur domain $\CO(H_i)$ of the handlebody $H_i$ is an open subset of $\CP\CM\CL(\D H_i)$, the space of projective measured laminations on $\D H_i$. If $\gamma$ is a weighted multicurve in the Masur domain then the pared manifold $(H_i,\supp(\gamma))$ has incompressible and acylindrical boundary, where $\supp(\gamma)$ is the support of $\gamma$ (see \cite{Masur86,Otal88} for the properties of the Masur domain). Kerckhoff \cite{Kerckhoff} proved that the Masur domain has full measure with respect to the measure class induced by the PL--structure of $\CP\CM\CL(\D H_i)$. The map $\phi:\D H_1\to \D H_2$ induces a homeomorphism $\phi_*:\CP\CM\CL(\D H_1)\to\CP\CM\CL(\D H_2)$ which preserves the canonical measure class. In particular, the intersection of $\CO(H_1)$ and $\phi_*^{-1}(\CO(H_2))$ is not empty and open in $\CP\CM\CL(\D H_1)$. Since weighted multicurves are dense in $\CP\CM\CL(\D H_1)$ the result follows.
\end{proof}

Now, choose a pant decomposition $P=\{p_1,\dots,p_{3g-3}\}$ of $\D H_1$ as in Lemma \ref{lem:OO} and identify it with a pant decomposition $P$ of $S$. Let $S\times[-2,2]$ be a regular neighborhood of $S$ in $M$ and $\CU$ a regular neighborhood of $P\times\{-1,1\}$ in $S\times[-2,2]$; $\CU$ is a union of disjoint open solid tori $U_1^+,\dots,U_{3g-3}^+,U_1^-,\dots,U_{3g-3}^-$ with $p_j^\pm=p_j\times\{\pm 1\}\subset U_j^\pm$ for all $j$. The boundary of the manifold $M^0=M\setminus\CU$ is a collection of tori
$$\D M^0=T_1^+\cup \dots\cup T_{3g-3}^+\cup T_1^-\cup\dots\cup T_{3g-3}^-$$ 
where $T_j^\pm$ bounds $U_j^\pm$. We choose a basis $(l_j^\pm,m_j^\pm)$ of $H_1(T_j^\pm;\BZ)$ for $j=1,\dots,3g-3$ as follows:

\begin{itemize}
\item[$l_j^\pm$:] For all $j$ there is a properly embedded annulus
$$\CA_j:(\BS^1\times[-1,1],\BS^1\times\{\pm 1\})\to(M^0\cap S\times [-2,2],T_j^\pm);$$
set $l_j^\pm=\CA_j\vert_{\BS^1\times\{\pm 1\}}$.
\item[$m_j^\pm$:] The curve $m_j^\pm$ is the meridian of the solid torus $U_j^\pm$ with the orientation chosen such that the algebraic intersection number of $m_j^\pm$ and $l_j^\pm$ is equal to 1.
\end{itemize}

For $n\in\BZ$ let $M^0_n$ be the manifold
$$M^0_n\defin M^0_{(nl_1^++m_1^+),\dots,(nl_{3g-3}^++m_{3g-3}^+),(-nl_1^-+m_1^-), \dots, (-nl_{3g-3}^-+m_{3g-3}^-)}$$
obtained by filling the curve $\pm nl_j^\pm+m_j^\pm$ for all $j$.

Let $V_j$ be a regular neighborhood of the image of $\CA_j$ in$M^0$; we may assume that $V_i\cap V_j=\emptyset$ for all $i\neq j$. The interior of the manifold $M^0\setminus\cup_jV_j$ is homeomorphic to $M\setminus P$ and its boundary is a collection $T_1,\dots,T_{3g-3}$ of tori. The complement of $M^0\setminus\cup_j V_j$ in $M^0_n$ is a union of $3g-3$ solid tori whose meridians do not depend on $n$. In particular, $M^0_n$ is homeomorphic to $M^0_0$ for all $n$. Since $M^0_0$ is, by construction, homeomorphic to $M$, we obtain

\begin{lem}\label{dehn}
The manifold $M^0_n$ is homeomorphic to $M$ for all $n\in\BZ$.\qed
\end{lem}

In order to complete the proof of Theorem \ref{main} we make use of the following result which we will show later on.

\begin{prop}\label{prop1}
There is a link $L\subset\inter M^0$ such that for all $\epsilon>0$ the manifold $\inter M^0$ admits a complete, finite volume hyperbolic cone--manifold structure with singular locus $L$ and cone--angle $2\pi-\epsilon$.
\end{prop}

We continue with the proof of Theorem \ref{main}. Since the manifold $\inter M^0$ admits a complete finite volume hyperbolic cone--manifold structure with cone--angle $2\pi-\epsilon$ it follows from the Generalized Dehn filling theorem that there is some $n$ such that $M^0_n$ admits a hyperbolic cone--manifold structure with cone--angle $2\pi-\epsilon$, too. This concludes the proof of  Theorem \ref{main} since $M$ and $M^0_n$ are homeomorphic by Lemma \ref{dehn}. 

We now prove Proposition \ref{prop1}. The surface $S$ separates $M^0$ in two manifolds $M^0_-$ and $M^0_+$. The boundary $\D M^0_\pm$ is the union of a copy of $S$ and the collection $\CP_\pm=\cup_{j=1,\dots,3g-3}T_j^{\pm 1}$ of tori. It follows from the choice of $P$ that the manifold $M^0_\pm$ is irreducible, atoroidal and has incompressible boundary. In particular, Thurston's Hyperbolization theorem \cite{Otal-Haken} implies that the interior of $M^0_\pm$ admits a complete hyperbolic metric with parabolic locus $\CP_\pm$. 

If $L\subset S$ is a simple closed curve such that $P\cup L$ fills $S$, then the pared manifold $(M^0_\pm,L)$ is acylindrical. Bonahon and Otal's Realization theorem implies that for all $\epsilon>0$ there is a geometrically finite hyperbolic metric $g_\pm$ on the interior of $M^0_\pm$ with parabolic locus $\CP_\pm$ and with bending lamination $\epsilon/2\cdot L$. The convex--core $CC(M^0_\pm,g_\pm)$ can be identified with $M^0_\pm\setminus\CP_\pm$ and hence the boundary of the convex--core consists of  a copy $S_\pm$ of the surface $S$; the identification of $S_\pm$ with $S$ induces a map $\psi:S_-\to S_+$ with
$$\inter M^0=CC(M^0_-,g_-)\cup_\psi CC(M^0_+,g_+).$$
The hyperbolic surface $S_\pm$ is bent along $L$ with dihedral angle $\frac 12\epsilon$. The following lemma concludes the proof of Proposition \ref{prop1}.

\begin{lem}\label{lem1}
The map $\psi:S_-\to S_+$ is isotopic to an isometry.
\end{lem}
\begin{proof}
The cover $(N_\pm,h_\pm)$ of $(\inter M^0_\pm,g_\pm)$ corresponding to the surface $S_\pm$ is geometrically finite. Since $S_\pm$ is incompressible we obtain that $N_\pm$ is homeomorphic to the interior of $S_\pm\times[-1,1]$ and the parabolic locus of $(N_\pm,h_\pm)$ is the collection $P\times\{\pm 1\}$. The convex surface $S_\pm\subset\inter M^0_\pm$ lifts to one of the components of the boundary of the convex--core of$(N_\pm,h_\pm)$; the other components are spheres with three punctures, and hence totally geodesic. The map $\psi$ can be extended to the map $\tilde\psi:N_-=S_-\times(-1,1)\to N_+=S_+\times(-1,1)$ given by $(x,t)\mapsto (\psi(x),-t)$. The map $\tilde\psi$ maps, up to isotopy, $\CP_-$ to $\CP_+$ and $L$ to $L$. Hence, the uniqueness part of Bonahon and Otal's Realization theorem implies that $\tilde\psi$ is isotopic to an isometry and this gives the desired result. 
\end{proof}

\sh{Concluding remarks}

Recall that in Theorem \ref{main} we do not claim that the singular locus of $d_\epsilon$ is independent of $\epsilon$. If $M$ is the double of a compact manifold with incompressible boundary whose interior admits a convex--cocompact hyperbolic metric, then, using Otal's trick, it is possible to construct a link $L$ such that $M$ admits a hyperbolic cone--manifold structure with singular locus $L$ and cone--angle $2\pi-\epsilon$ for all $\epsilon$. Proposition \ref{prop1} suggests that this may be a more general phenomenon but the author does not think that it is always possible to choose the singular locus independently of $\epsilon$.

\medskip

{\bf Question 1}\qua {\sl Let $L$ be a link in $\BS^2\times\BS^1$ which	intersects an essential sphere $n$ times. Is there	a hyperbolic cone--manifold structure on $\BS^2\times\BS^1$ with	singular locus $L$ and with cone--angle greater than $\frac{n-2}n2\pi$?}

\medskip

{\bf Question 2}\qua {\sl Is there a link $L\subset\BS^3$ such that	for every $\epsilon>0$ there is a hyperbolic cone--manifold structure	on $\BS^3$ with singular locus $L$ and with cone--angle $2\pi-\epsilon$?}

\medskip

We suspect that both questions have a negative answer.

We define, as suggested by Michel Boileau, the {\em hyperbolic volume} $\hypvol(M)$ of a closed 3--manifold $M$ as the infimum of the volumes of all possible hyperbolic cone--manifold structures on $M$ with cone--angle less or equal to $2\pi$. It follows from \cite{Hodgson-Kerckhoff} and from the Schl\"afli formula that the hyperbolic volume of a manifold $M$ is achieved if and only if $M$ is hyperbolic. A sequence of hyperbolic cone--manifold structures realizes the hyperbolic volume if the associated volumes converge to $\hypvol(M)$. From the arguments used in the proof of the Orbifold theorem \cite{Boileau-Porti} it is easy to deduce that the hyperbolic volume is realized by a sequence of hyperbolic cone--manifold structures whose cone--angles are all greater or equal to $\pi$.   
\medskip

{\bf Question 3}\qua {\sl Is there a sequence of metrics realizing the hyperbolic volume and such that the associated cone--angles tend to $2\pi$?}

\medskip

As remarked in the introduction, it follows from our construction that there are sequences of hyperbolic cone--manifold structures whose cone--angles tend to $2\pi$ and which have uniformly bounded volume.

Let $M$ now be a closed orientable and irreducible 3--manifold $M$. We say that $M$ is {\em geometrizable} if Thurston's Geometrization Conjecture holds for it. If $M$ is geometrizable then let $M_{\hyp}$ be the associated complete finite volume hyperbolic manifold. In \cite{Juan-thesis} we proved:

\begin{sat*}
Let $M$ be a closed, orientable, geometrizable and prime 3--manifold. Then the minimal volume $\Minvol(M)$ of $M$ is equal to $\vol(M_{\hyp})$ and moreover, the manifolds $(M,g_i)$ converge in geometrically to $M_{\hyp}$ for every sequence $(g_i)$ of metrics realizing $\Minvol(M)$. In particular, the minimal volume is achieved if and only if $M$ is hyperbolic.
\end{sat*}

Recall that the minimal volume $\Minvol(M)$ of $M$ is the infimum of the volumes $\vol(M,g)$ of all Riemannian metrics $g$ on $M$ with sectional curvature bounded in absolute value by one. A sequence of metrics $(g_i)$ realizes the minimal volume if their sectional curvatures are bounded in absolute value by one and if $\vol(M,g_i)$ converges to $\Minvol(M)$.

Under the assumption that the manifold $M$ is geometrizable and prime, it follows with the same arguments as in \cite{Juan-thesis} that the hyperbolic volume can be bounded from below by the minimal volume.
\medskip

{\bf Question 4}\qua {\sl If $M$ is geometrizable and prime, do the hyperbolic and the minimal volume coincide?}

\medskip

This question has a positive answer if the manifold $M$ is the double of a manifold which admits a convex--cocompact metric and the answer should be also positive without this restriction. If this is the case, then it should also be possible to show that the Gromov--Hausdorff limit of every sequence of hyperbolic cone--manifold structures which realizes the hyperbolic volume is isometric to $M_{\hyp}$. We do not dare to ask if the assumption on $M$ to be geometrizable can be dropped.

\bibliographystyle{gtart}

\end{document}